\documentclass[11pt]{article}
\usepackage{latexsym,amsfonts,amsthm,amsmath,amscd,amssymb,color,url}
\usepackage[dvips]{graphicx}

\newcommand\matP{{\mathbb{P}}}

\newcommand\matC{{\mathbb{C}}}

\newcommand\Dhat{{\widehat{D}}}

\newcommand\Sigmatil{{\widetilde\Sigma}}

\newcommand\Btil{{\widetilde B}}
\newcommand\Rtil{{\widetilde R}}

\renewcommand{\hbar}{{\overline{h}}}

\newfont{\Got}{eufm10 scaled 1200}

\newcommand{\mycap} [1] {\caption{\footnotesize{#1}}}

\begin{document}

\title{Elementary solution of an infinite sequence of instances of the Hurwitz problem}

\author{Tom~\textsc{Ferragut}\and Carlo~\textsc{Petronio}}

\maketitle

\begin{abstract}\noindent
We prove that there exists no branched cover from the torus to the sphere with degree $3h$
and $3$ branching points in the target with local degrees $(3,\ldots,3),\ (3,\ldots,3),\
(4,2,3,\ldots,3)$ at their preimages. The result was already established
by Izmestiev, Kusner, Rote, Springborn, and Sullivan,
using geometric techniques, and by Corvaja and Zannier with a more algebraic approach,
whereas our proof is topological and completely elementary: besides the definitions, it
only uses the fact that on the torus a simple closed curve can only be \emph{trivial}
(in homology, or equivalently bounding a disc, or equivalently separating) or \emph{non-trivial}.
\smallskip

\noindent MSC (2010): 57M12.
\end{abstract}

\noindent
A (topological) branched cover between surfaces is a
map $f:\Sigmatil\to\Sigma$, where $\Sigmatil$
and $\Sigma$ are closed and connected $2$-manifolds and $f$ is locally modeled
(in a topological sense) on maps of the form
$(\matC,0)\ni z\mapsto z^k\in(\matC,0)$. If $k>1$ the point
$0$ in the target $\matC$ is called a \emph{branching point},
and $k$ is called the local degree at the point $0$ in the source $\matC$.
There are finitely many branching points, removing which, together
with their pre-images, one gets a genuine cover of some degree $d$.
If there are $n$ branching points, the local degrees at the points
in the pre-image of the $j$-th one form a partition $\pi_j$ of $d$ of some
length $\ell_j$, and the following Riemann-Hurwitz relation holds:
$$\chi\left(\Sigmatil\right)-(\ell_1+\ldots+\ell_n)=d\left(\chi\left(\Sigma\right)-n\right).$$
The very old \emph{Hurwitz problem} asks whether given $\Sigmatil,\Sigma,d,n,\pi_1,\ldots,\pi_n$
satisfying this relation there exists some $f$ realizing them. (For a
non-orientable $\Sigmatil$ and/or $\Sigma$ the Riemann-Hurwitz relation
must actually be complemented with certain other necessary conditions,
but we will not get into this here.) A number of partial solutions
of the Hurwitz problem have been obtained over the time, and we quickly
mention here the fundamental~\cite{EKS}, the survey~\cite{Bologna}, and
the more recent~\cite{Pako, PaPe, CoPeZa, PaPebis, SongXu}.

\smallskip

Certain instances of the Hurwitz problem recently emerged in
the work of M.~Zieve~\cite{Zieve} and his team of collaborators, including
in particular the case
where the source surface is the torus $T^2$, the
target is the sphere $S^2$, the degree is $d=3h$,
and there are $n=3$ branching points with associated partitions
$(3,\ldots,3),\ (3,\ldots,3),\ (4,2,3,\ldots,3)$ of $d$.
It actually turns out that this branch datum is indeed not realizable, as Zieve
had conjectured, which follows from results established in~\cite{GeomDed}
using geometric techniques (holonomy of Euclidean structures).
The same fact was also elegantly proved by Corvaja and Zannier~\cite{CoZa} with a more
algebraic approach.  In this note we provide yet another proof of the same
result.  Our approach is purely combinatorial and completely elementary:
besides the definitions, it
only uses the fact that on the torus a simple closed curve can only be \emph{trivial}
(in homology, or equivalently bounding a disc, or equivalently separating) or \emph{non-trivial}.

\medskip

We conclude this introduction with the formal statement of the (previously known)
result established in this note:

\medskip

\noindent
\textbf{Theorem}. There exists no branched cover $f:T^2\to S^2$ with degree $d=3h$ and
$3$ branching points with associated partitions $$(3,\ldots,3),\ (3,\ldots,3),\
(4,2,3,\ldots,3).$$

\section{Dessins d'enfant}

In this section we quickly review the beautiful technique of
dessins d'enfant due to Grothendieck~\cite{Cohen, Groth},
noting that, at the elementary level at which we exploit it, it
only requires the definition of branched cover and some very
basic topology.

\smallskip

\begin{figure}
    \begin{center}
    \includegraphics[scale=.6]{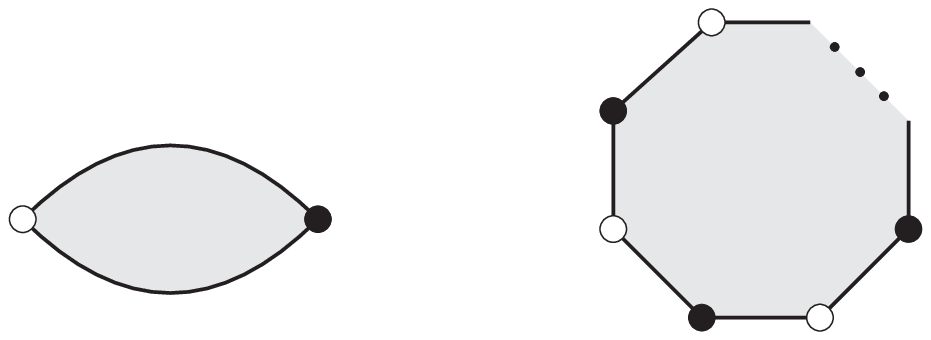}
    \end{center}
\mycap{A bigon and a polygon.\label{regions:fig}}
\end{figure}

Let $f:\Sigmatil\to S^2$ be a degree-$d$ branched cover from a closed connected surface $\Sigmatil$ to the sphere $S^2$, branched over $3$ points
$p_1,p_2,p_3$ with local degrees $\pi_j=\left(d_{ji}\right)_{i=1}^{\ell_j}$ over $p_j$.
In $S^2$ take a simple arc $\sigma$ with vertices at $p_1$ (white) and $p_2$ (black), and we view $S^2$ as being obtained from the
(closed) bigon $\Btil$ of Fig.~\ref{regions:fig}-left by attaching both the edges of $\Btil$ to $\sigma$ so to match the vertex colors.
This gives a realization of $S^2$ as the quotient of $\Btil$ under the identification of its two edges.
Let $\lambda:\Btil\to S^2$ be the projection to the quotient.
Note that the complement of $\sigma$ in $S^2$ is an open disc $B$,
whose closure in $S^2$ is the whole of $S^2$, but the restriction of $\lambda$ to the interior of $\Btil$ is a homeomorphism with $B$, so we can
view $\Btil$ as the abstract closure of $B$.

Now set $D=f^{-1}(\sigma)$. Then $D$ is a graph with white vertices of valences
$\left(d_{1i}\right)_{i=1}^{\ell_1}$ and black vertices of valences $\left(d_{2i}\right)_{i=1}^{\ell_2}$, and
$D$ is bipartite (every edge has a white and a black end). Moreover the complement of $D$ in $\Sigmatil$ is a union of open discs
$\left(R_i\right)_{i=1}^{\ell_3}$, where $R_i$ is the interior of a polygon with $2d_{3i}$ vertices of alternating white and black color.
This means that, if $\Rtil_i$ is the polygon of Fig.~\ref{regions:fig}-right (with $2d_{3i}$ vertices),
there exists a map $\lambda_i:\Rtil_i\to\Sigmatil$ which restricted to the interior of $\Rtil_i$ is a homeomorphism with $R_i$,
and restricted to each edge is a homeomorphism with an edge of $D$ matching the vertex colors. So $\Rtil_i$ can be viewed as
the abstract closure of $R_i$. The map $\lambda_i$ may fail to be a homeomorphism between $\Rtil_i$ and
the closure of $R_i$ in $\Sigmatil$ if $R_i$ is multiply incident
to some vertex of $D$ or doubly incident to some edge of $D$. We say that $R_i$ has embedded closure if $\lambda_i$ is injective,
hence a homeomorphism between $\Rtil_i$ and the closure of $R_i$ in $\Sigmatil$.

\medskip

We will way that a bipartite graph $D$ in $\Sigmatil$ with valences $\left(d_{1i}\right)_{i=1}^{\ell_1}$
at the white vertices and $\left(d_{2i}\right)_{i=1}^{\ell_2}$ at the black ones, and
complement consisting of polygons having $\left(2d_{3i}\right)_{i=1}^{\ell_3}$ edges,
\emph{realizes} the branched cover $f:\Sigmatil\to S^2$ with $3$ branching points and
local degrees $\pi_1,\pi_2,\pi_3$ over them. This terminology is justified
by the fact that $f$ exists if and only if $D$ does.

\section{Proof of the Theorem}

\begin{figure}
    \begin{center}
    \includegraphics[scale=.6]{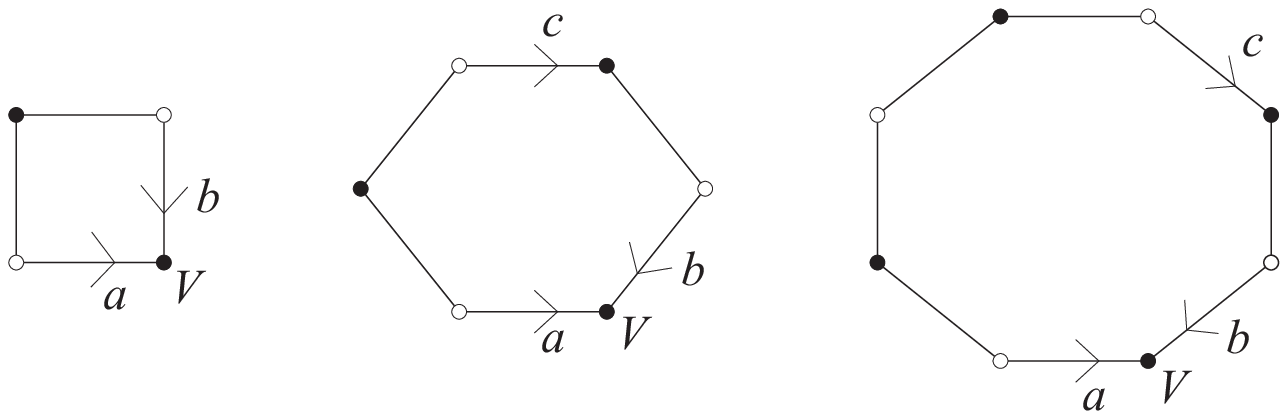}
    \end{center}
\mycap{The regions. The notation $V,a,b,c,$ is only needed for the base step of our induction argument.\label{SHO:fig}}
\end{figure}

Suppose by contradiction that a branched cover $f:T^2\to S^2$ as in the statement exists,
and let $D$ be a dessin d'enfant on  $T^2$ realizing it,
as explained in the previous section,
with white and black vertices corresponding to the first two partitions, so the complementary regions are one square $S$,
some hexagons $H$ and one octagon $O$, shown abstractly in Fig.~\ref{SHO:fig}.
Let $\Dhat$ be the graph dual to $D$
(which is well-defined because the complement of $D$ is a union of open discs),
and let $\Gamma$ be the set of all simple loops in $\Dhat$
which are simplicial (concatenations of edges), and non-trivial
(non-zero in $H_1(T^2)$, or, equivalently, not bounding a disc on $T^2$, or, equivalently, not separating $T^2$).
Since $T^2\setminus\Dhat$ is also a union of open discs, the inclusion $\Dhat\hookrightarrow T^2$
induces a surjection $H_1(\Dhat)\to H_1(T^2)$. Moreover $H_1(\Dhat)$ is generated by simple simplicial loops,
so $\Gamma$ is non-empty. We now define
$\Gamma_n$ as the set of loops in $\Gamma$ consisting of $n$ edges, and we prove by induction that $\Gamma_n=\emptyset$,
thereby showing that $\Gamma=\emptyset$ and getting the desired contradiction.

\begin{figure}
    \begin{center}
    \includegraphics[scale=.6]{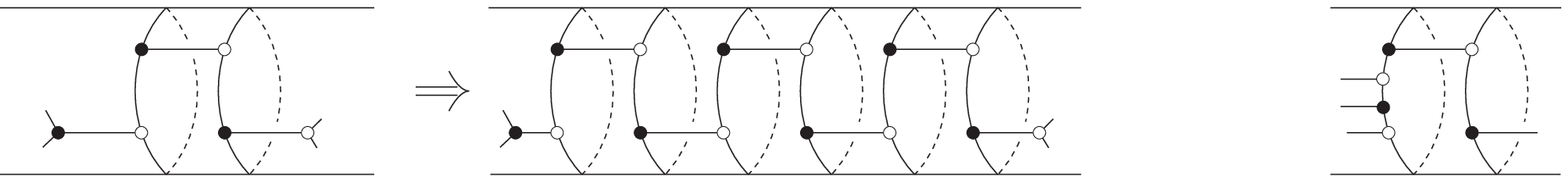}
    \end{center}
\mycap{A tube of $H$'s cannot end at an $O$ from both sides. A non-embedded $O$ gives a non-embedded $H$\label{Htube:fig}}
\end{figure}

For $n=1$ we prove the slightly stronger fact (needed below) that every region has embedded closure, namely,
that its closure in $T^2$ is homeomorphic to its abstract closure.
Taking into account the symmetries (including a color switch) this may fail to happen
only if some edge $a$ in Fig.~\ref{SHO:fig} is glued to $b$ or $c$ of the same region
(if two vertices of a region are glued together then two edges also are, since the vertices have valence 3).
The case $b=a$ implies $V$ has valence $1$, so it is impossible. If $c=a$ in $H$ we have the
situation of Fig.~\ref{Htube:fig}-left, and each of the neighboring regions
already has $3$ vertices of one color, so it cannot be $S$. If it is an $H$, it
also has a gluing of type $c=a$. Iterating, we have a tube of $H$'s as in Fig.~\ref{Htube:fig}-centre that
at some point must hit $O$ from both sides, which is impossible because the terminal region already contains
$5$ vertices of each color. If $c=a$ in $O$ then we have Fig.~\ref{Htube:fig}-right, so a neighboring
region also has non-embedded closure, which was already excluded.

\begin{figure}
    \begin{center}
    \includegraphics[scale=.6]{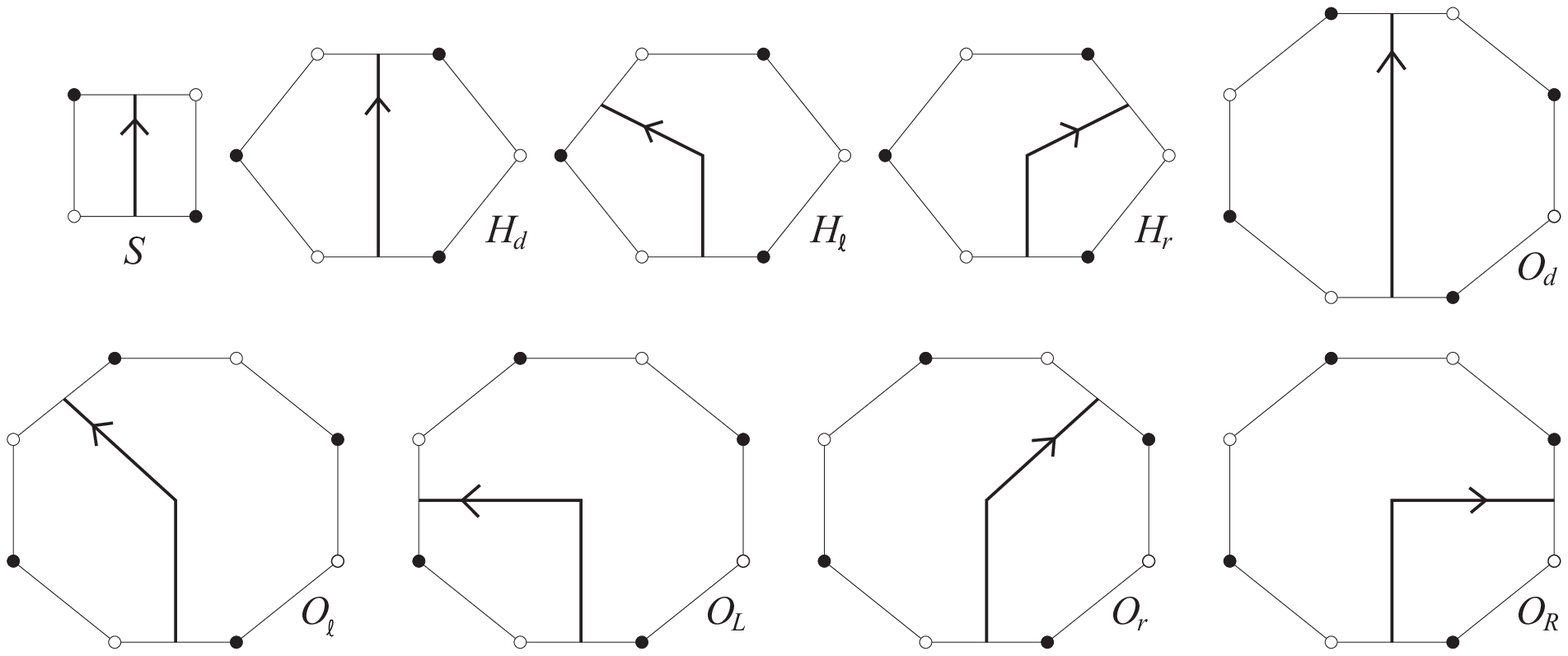}
    \end{center}
\mycap{The ways $\gamma$ can cross a region. Note that the vertex colors may be switched. \label{Crossings:fig}}
\end{figure}

Let us now assume that $n\geqslant 2$ and $\Gamma_m=\emptyset$ for all $m<n$.
By contradiction, take $\gamma\in\Gamma_n$. From now on in our figures we will use for $\gamma$ a thicker line than
that used for $D$. We first note that $\gamma$ cannot enter a region through an edge and leave it from an
adjacent edge (otherwise we could reduce its length), so the only ways $\gamma$ can cross a region
are those shown in Fig.~\ref{Crossings:fig}.
Therefore $\gamma$ is described by a word in the letters $S,H_d,H_\ell,H_r,O_d,O_\ell,O_L,O_r,O_R$, from which
we omit the $H_d$'s for simplicity. The vertex coloring implies that the total number of $S,H_\ell,H_r,O_d,O_L,O_R$ in $\gamma$ is even.

\begin{figure}
    \begin{center}
    \includegraphics[scale=.33]{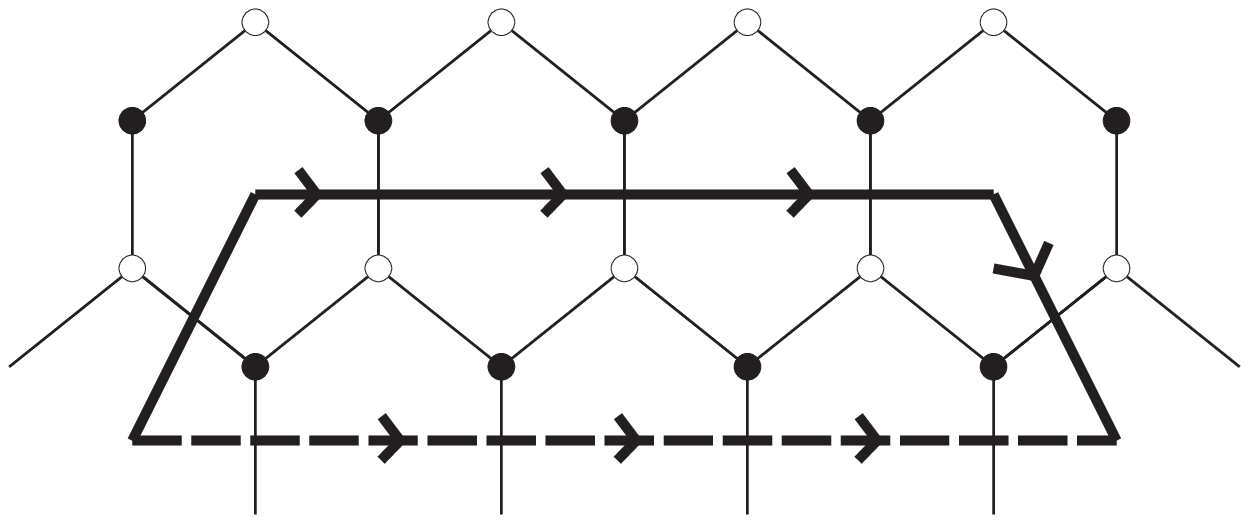}\quad\quad
    \includegraphics[scale=.33]{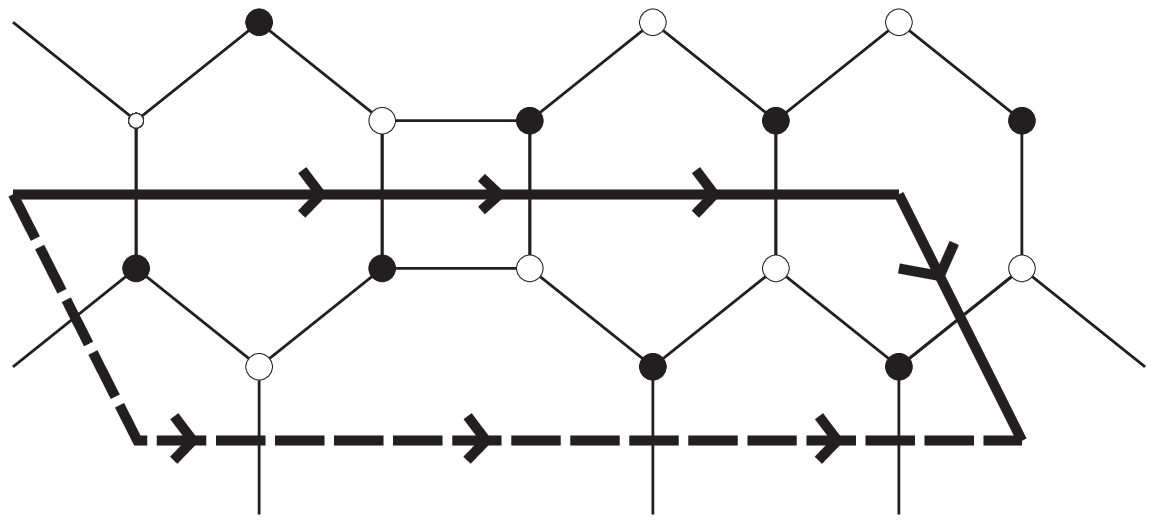}
    \end{center}
\mycap{Impossible configurations. \label{ImpS:fig}}
\end{figure}

We now prove that any subword $H_rH_r$, $H_\ell H_\ell$, $SH_r$ or $SH_\ell$ is impossible in $\gamma$, as shown in
Fig.~\ref{ImpS:fig}
(here the thick dashed line gives a new $\gamma$ contradicting the minimality of the original one).
This already implies the former of the following claims:

\begin{itemize}
\item[$(1)$] No $\gamma\in\Gamma_n$ can contain $S$ but not $O$;
\item[$(2)$] There exists $\gamma\in\Gamma_n$ consisting of $H$'s only.
\end{itemize}

To establish the latter, we suppose $O\in\gamma\in\Gamma_n$ and list all the possible cases up to symmetry
(which includes switching colors and/or reversing the direction of $\gamma$):

\begin{center}
\begin{tabular}{l|l|l}
    $S\notin\gamma$ & $O_d\in\gamma\Rightarrow\gamma=O_dH_r(H_\ell H_r)^p$ & \\ \cline{2-3}
                    & $O_r\in\gamma\Rightarrow\gamma=O_r(H_r H_\ell)^p$    &  \\ \cline{2-3}
                    & $O_R\in\gamma$                                   & $\gamma=O_RH_r(H_\ell H_r)^p$ \\ \cline{3-3}
                    &                                                  & $\gamma=O_RH_\ell(H_rH_\ell)^p$ \\ \hline
    $S\in\gamma\Rightarrow SO\in\gamma$
                      & $O_d\in\gamma\Rightarrow\gamma=SO_d(H_rH_\ell)^p$ & \\ \cline{2-3}
                      & $O_r\in\gamma$                                   & $\gamma=H_rSO_r(H_rH_\ell)^p$ \\ \cline{3-3}
                      &                                                  & $\gamma=H_\ell SO_r(H_\ell H_r)^p$ \\ \cline{2-3}
                      & $O_R\in\gamma$                                   &
\end{tabular}
\end{center}

\begin{figure}
    \begin{center}
    \includegraphics[scale=.33]{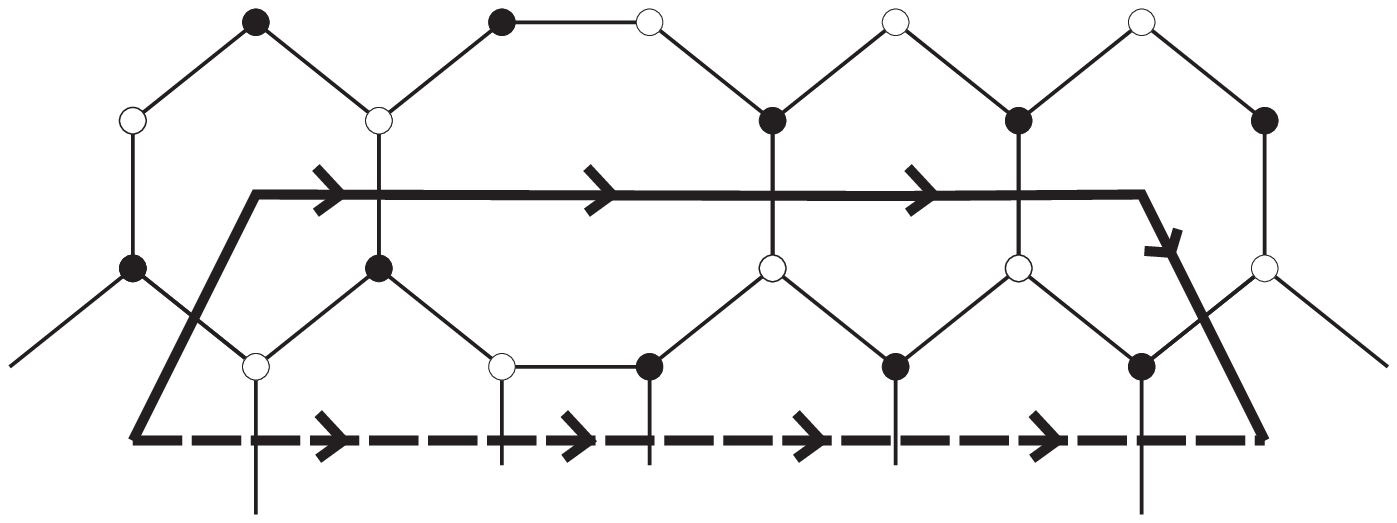}\quad
    \includegraphics[scale=.33]{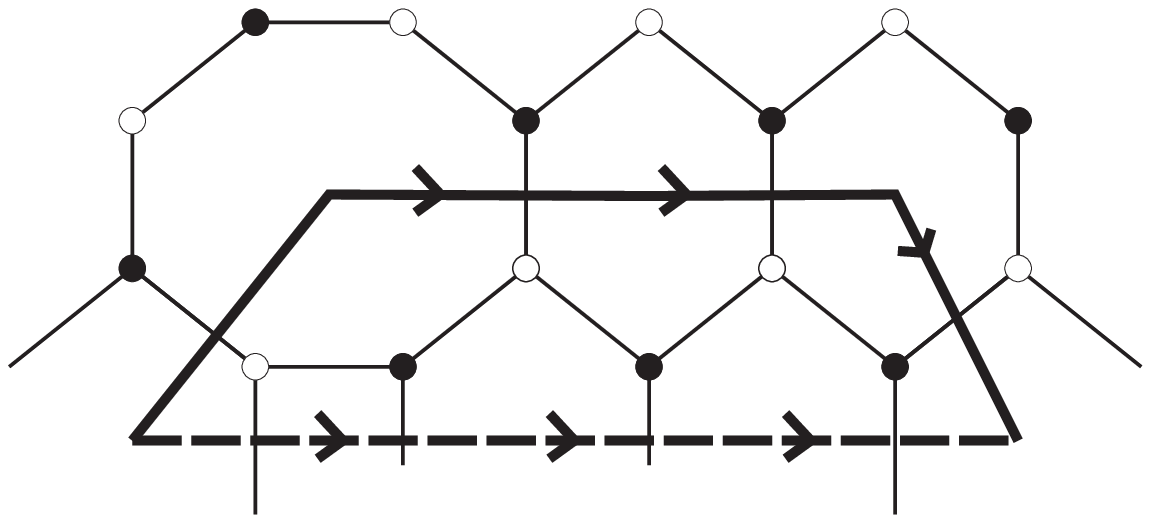}\quad
    \includegraphics[scale=.33]{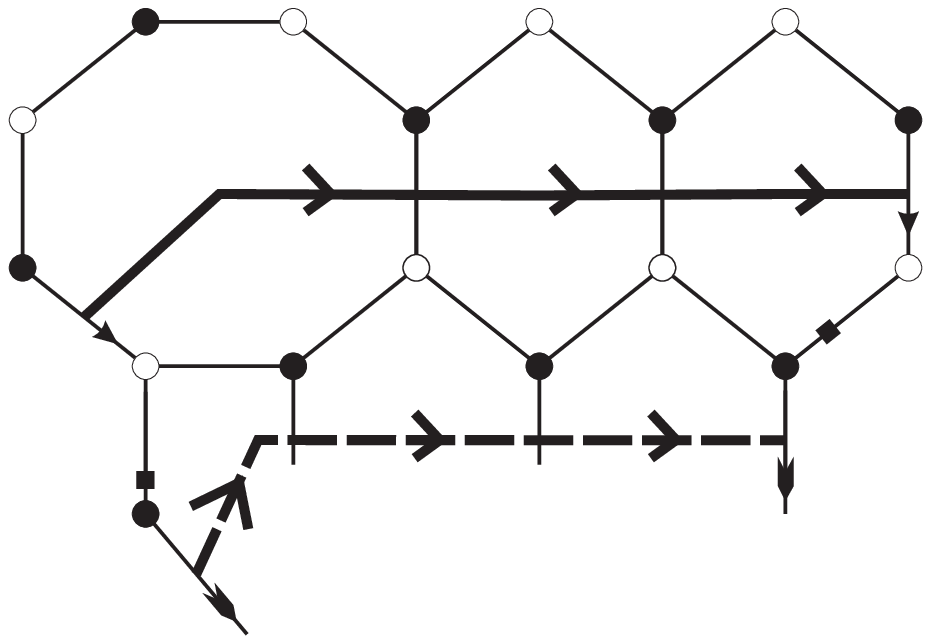}
    \end{center}
\mycap{$\gamma=O_dH_r(H_\ell H_r)^p\Rightarrow\exists\gamma'\ldots$ (left);
$\gamma=O_r(H_r H_\ell)^p\Rightarrow\exists\gamma'\ldots$ (centre for $p>0$ and right for $p=0$).
On the right, as in many figures below, we decorate some edges to indicate that they are glued in pairs.\label{1ab:fig}}
\end{figure}

\begin{figure}
    \begin{center}
    \includegraphics[scale=.33]{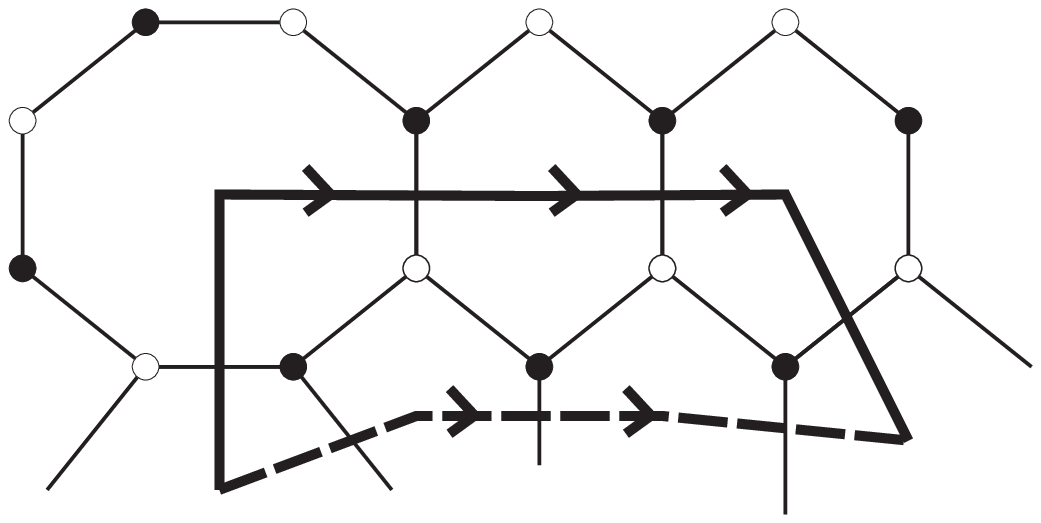}\quad\quad
    \includegraphics[scale=.33]{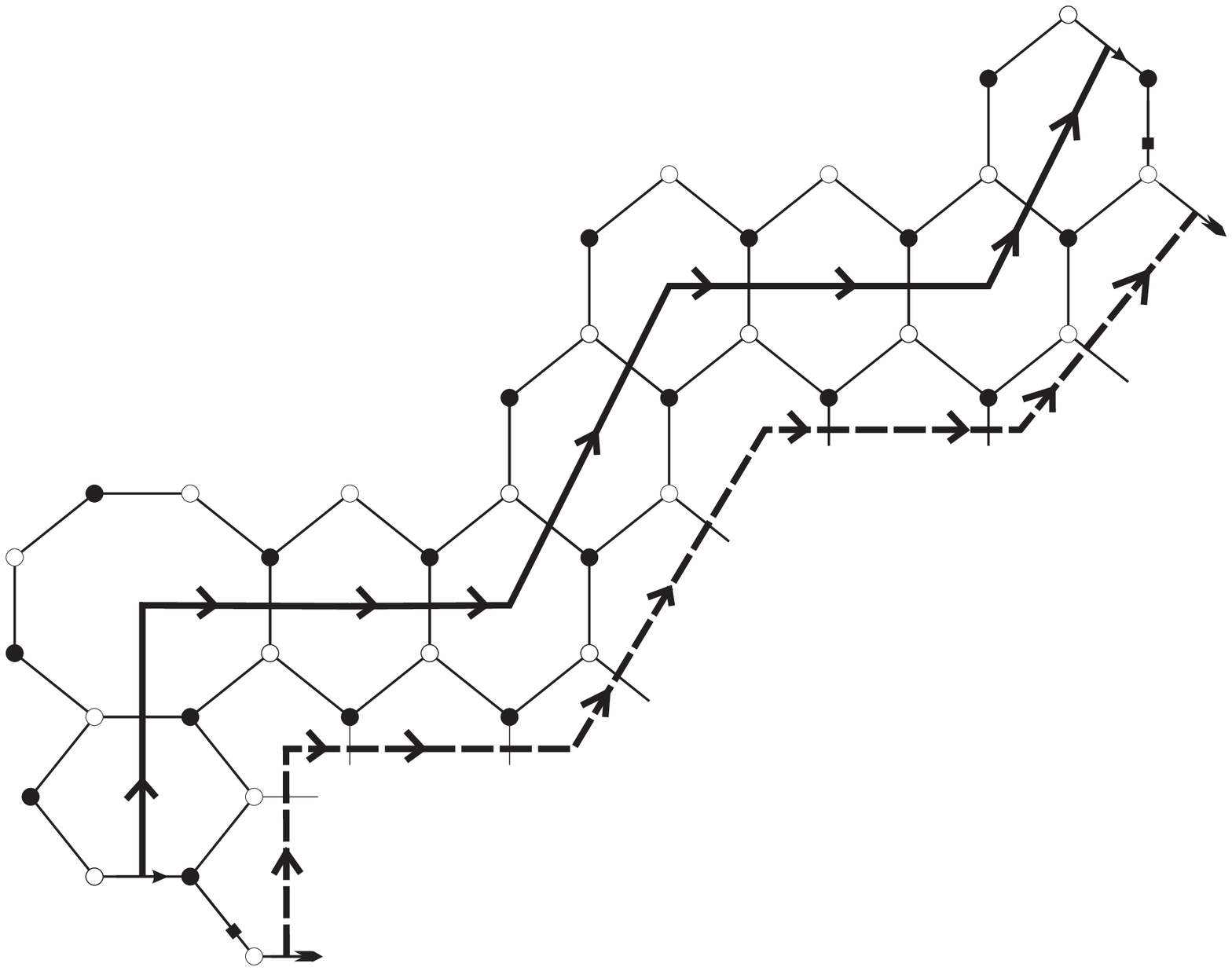}
    \end{center}
\mycap{$\gamma=O_RH_r(H_\ell H_r)^p$ impossible (left); $\gamma=O_RH_\ell(H_rH_\ell)^p\Rightarrow\exists\gamma'\ldots$ (right). \label{1c:fig}}
\end{figure}

\begin{figure}
    \begin{center}
    \includegraphics[scale=.33]{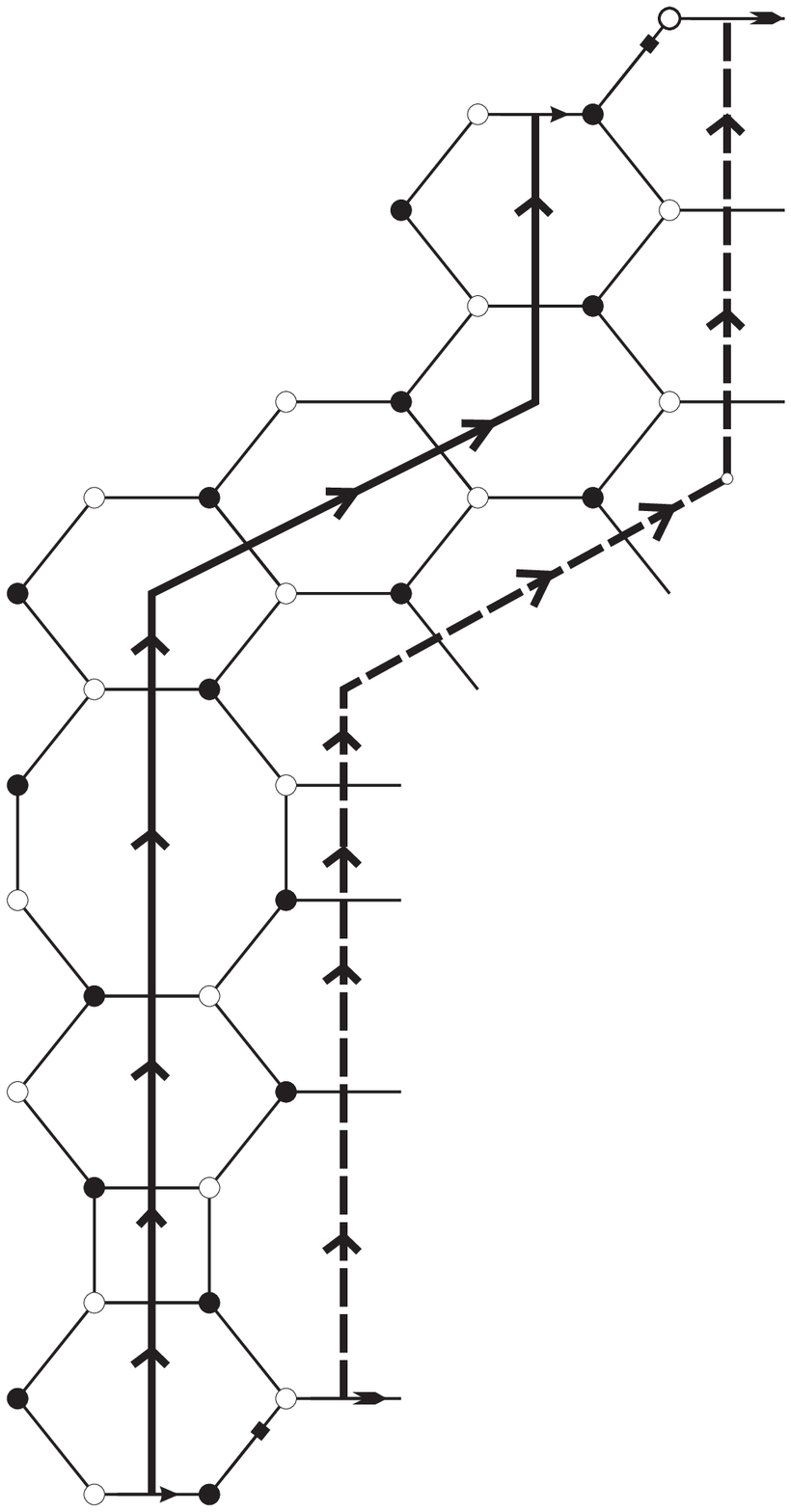}\quad\quad
    \includegraphics[scale=.33]{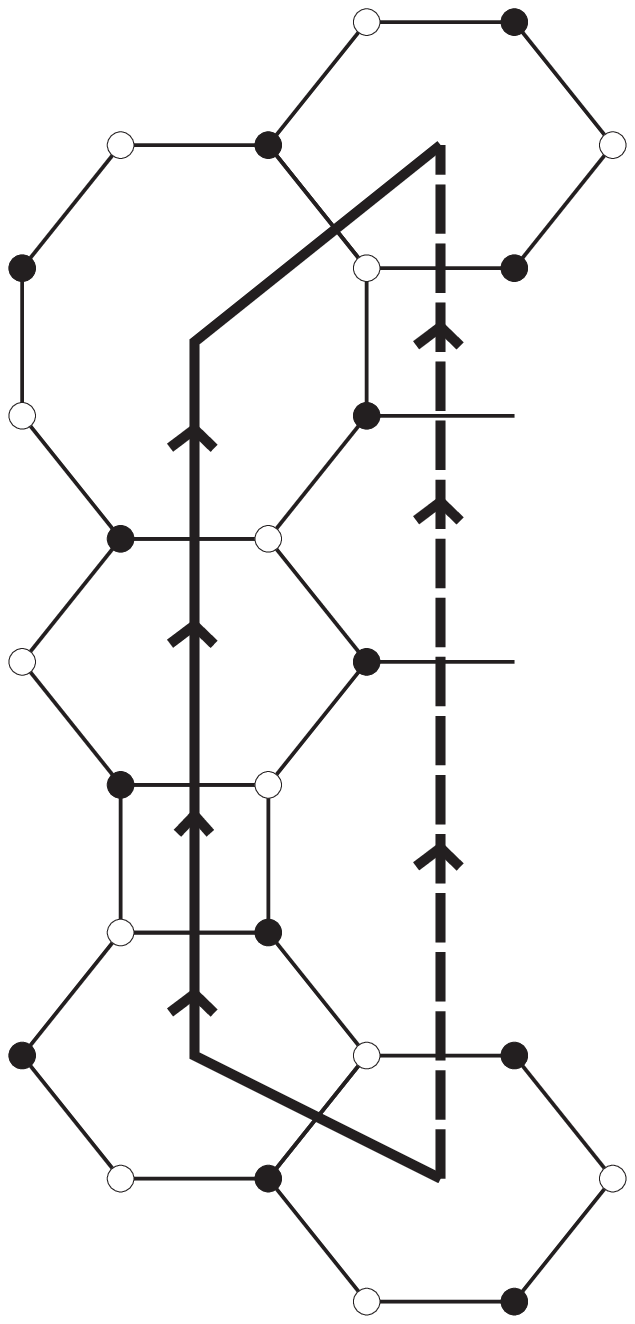}\quad\quad
    \includegraphics[scale=.33]{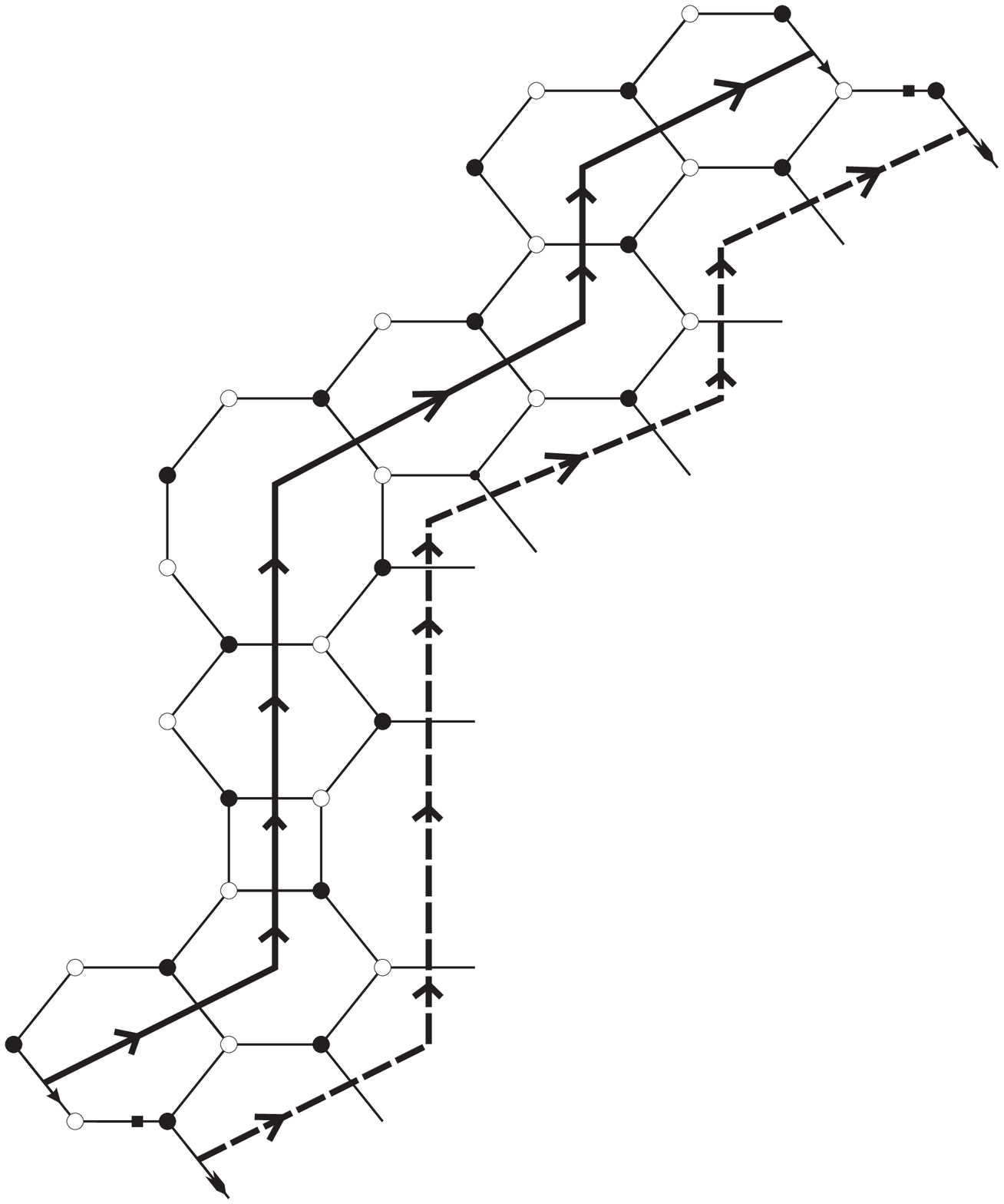}
    \end{center}
\mycap{$\gamma=SO_d(H_rH_\ell)^p\Rightarrow\exists\gamma'\ldots$ (left);
$\gamma=H_rSO_r(H_rH_\ell)^p$ impossible (centre);
$\gamma=H_\ell SO_r(H_\ell H_r)^p\Rightarrow\exists\gamma'\ldots$ (right). \label{2ab:fig}}
\end{figure}

\begin{figure}
    \begin{center}
    \includegraphics[scale=.33]{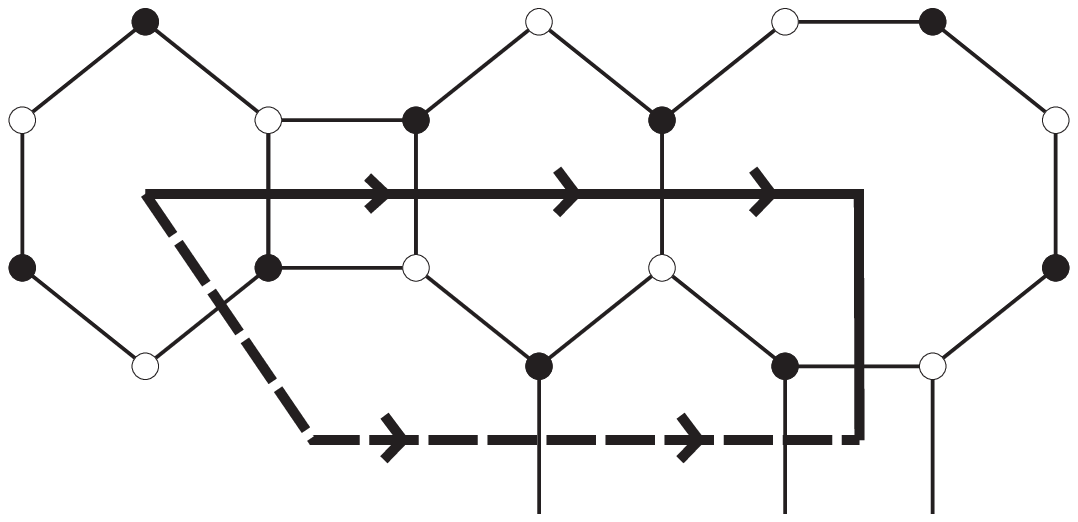}\quad\quad
    \includegraphics[scale=.33]{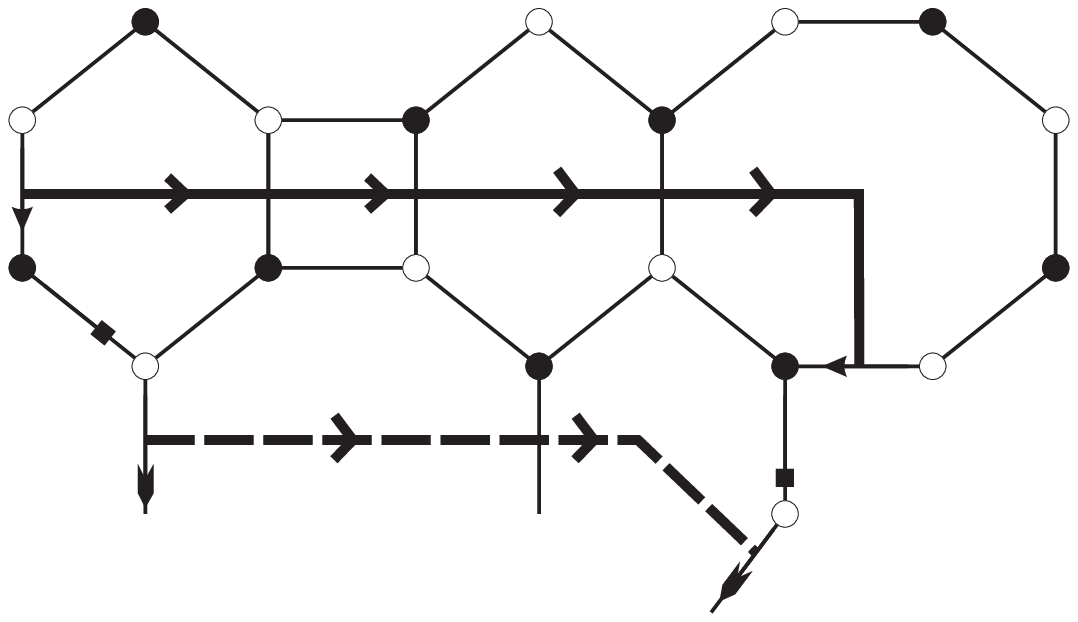}\quad\quad
    \includegraphics[scale=.33]{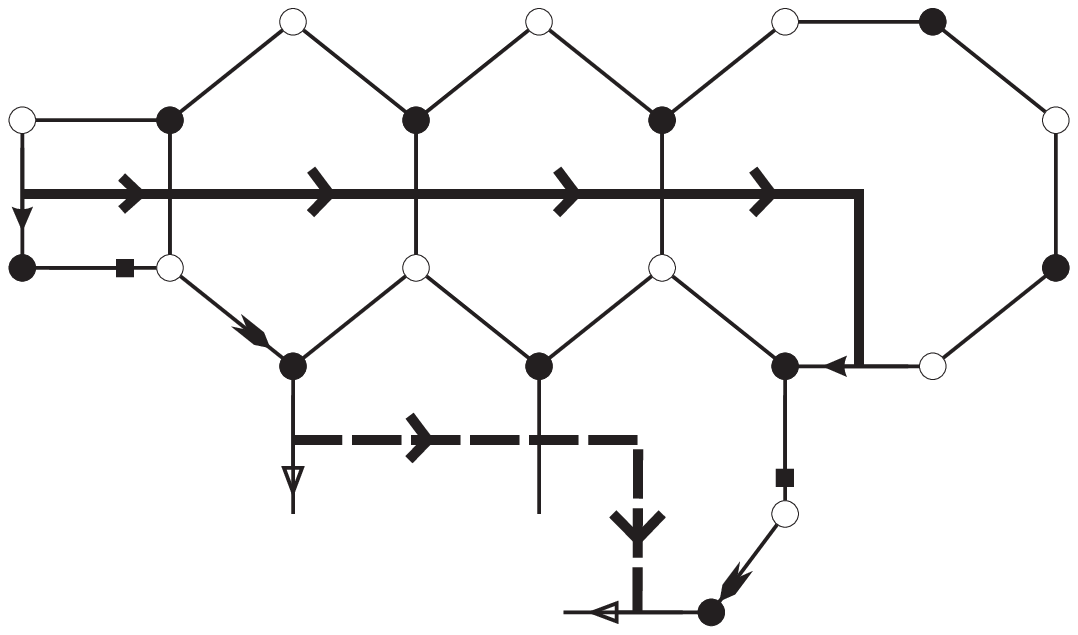}
    \end{center}
\mycap{$SO_R\in\gamma$ impossible: if in $\gamma$ there are $m$ copies of $H$ outside the word $SO_R$,
we treat separately the cases $m\geqslant 2$ (left), $m=1$ (centre), and $m=0$ (right).
In the last case, if there are not even $H_d$'s between $S$ and $O_R$, the absurd comes from the fact that a bigon is created.
\label{2c:fig}}
\end{figure}

\begin{figure}
    \begin{center}
    \includegraphics[scale=.33]{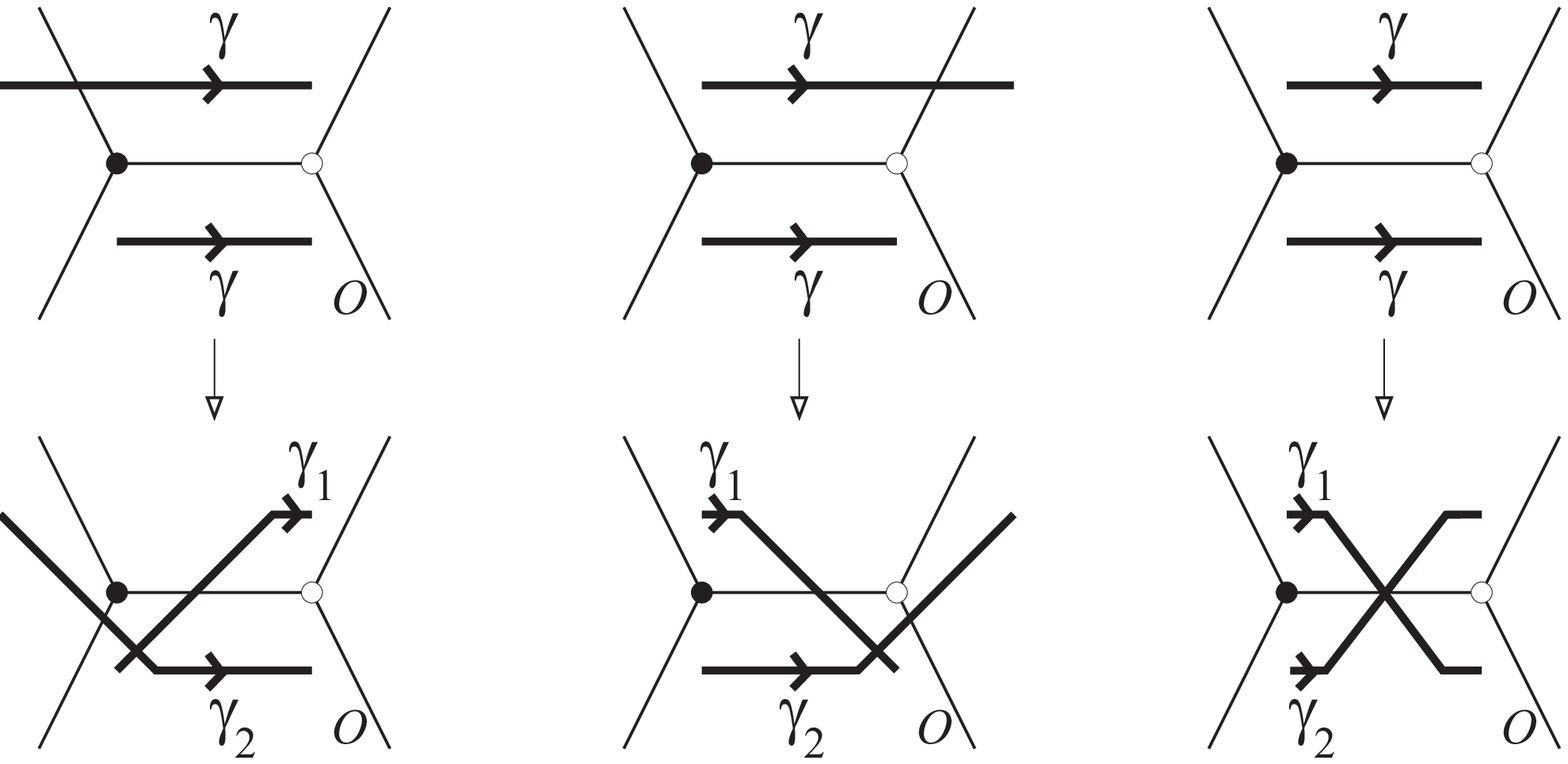}
    \end{center}
\mycap{If $O$ appears in $\gamma$ and immediately to the right of $\gamma$, we can construct two loops
$\gamma_{1,2}$ of which $\gamma$ is the homological sum.\label{gamma12moves:fig}}
\end{figure}

For each of these cases we show in Figg.~\ref{1ab:fig} to~\ref{2c:fig} a modification of $\gamma$
which gives a new loop $\gamma'$ isotopic to $\gamma$ (and hence in $\Gamma$), and not longer
than $\gamma$. When $\gamma'$ is shorter than $\gamma$ we have a contradiction to the minimality of $\gamma$, so the
case is impossible. To conclude we must show that $\gamma'$ does not contain $O$
in the cases where it is as long as $\gamma$. To do this, suppose that $\gamma'$ contains $O$,
and construct two loops $\gamma_{1,2}$ by applying
one of the three moves of Fig.~\ref{gamma12moves:fig}.
Note that whatever move applies, $\gamma$ is the homological sum of $\gamma_1$ and $\gamma_2$, so at least
one of them is non-trivial. If one of the moves of Fig.~\ref{gamma12moves:fig}-left/centre applies,
the total length of $\gamma_1$ and $\gamma_2$ is $1$ plus the length of $\gamma$, but
we know that there is no length-1 loop at all (trivial or not), so both $\gamma_1$ and $\gamma_2$
are shorter than $\gamma$, a contradiction. If only the move of Fig.~\ref{gamma12moves:fig}-right applies
then we are either in Fig.~\ref{1c:fig}-right, or Fig.~\ref{2ab:fig}-left or Fig.~\ref{2ab:fig}-right
and $O$ is the region where $\gamma'$ makes a left turn; in this case
the total length of $\gamma_1$ and $\gamma_2$ is $2$ plus the length of $\gamma$, but
$\gamma_1$ and $\gamma_2$ both have length at least $3$, so they are both shorter than $\gamma$, and
again we have a contradiction.

\begin{figure}
    \begin{center}
    \includegraphics[scale=.33]{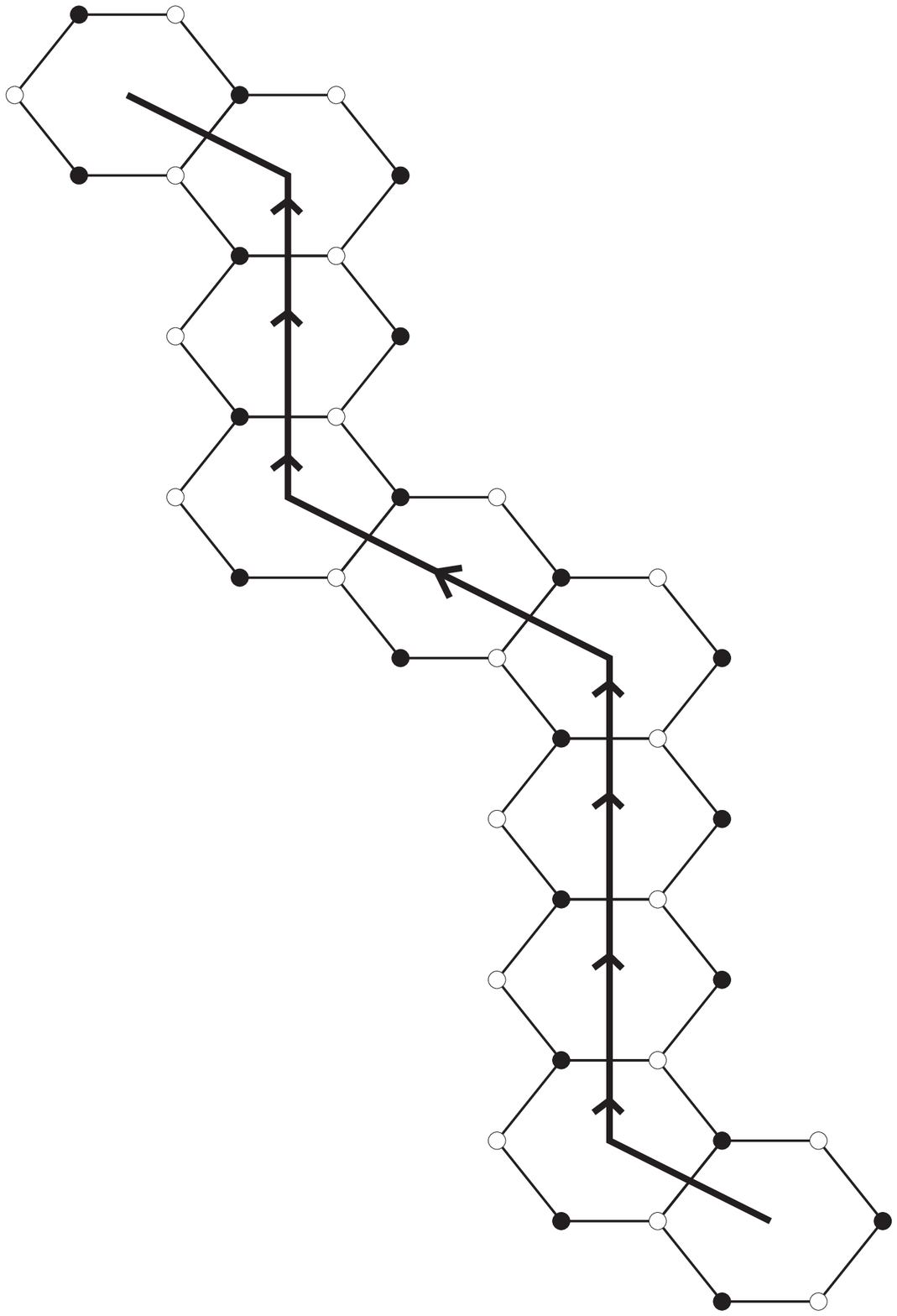}\hspace{-.5cm}
    \includegraphics[scale=.33]{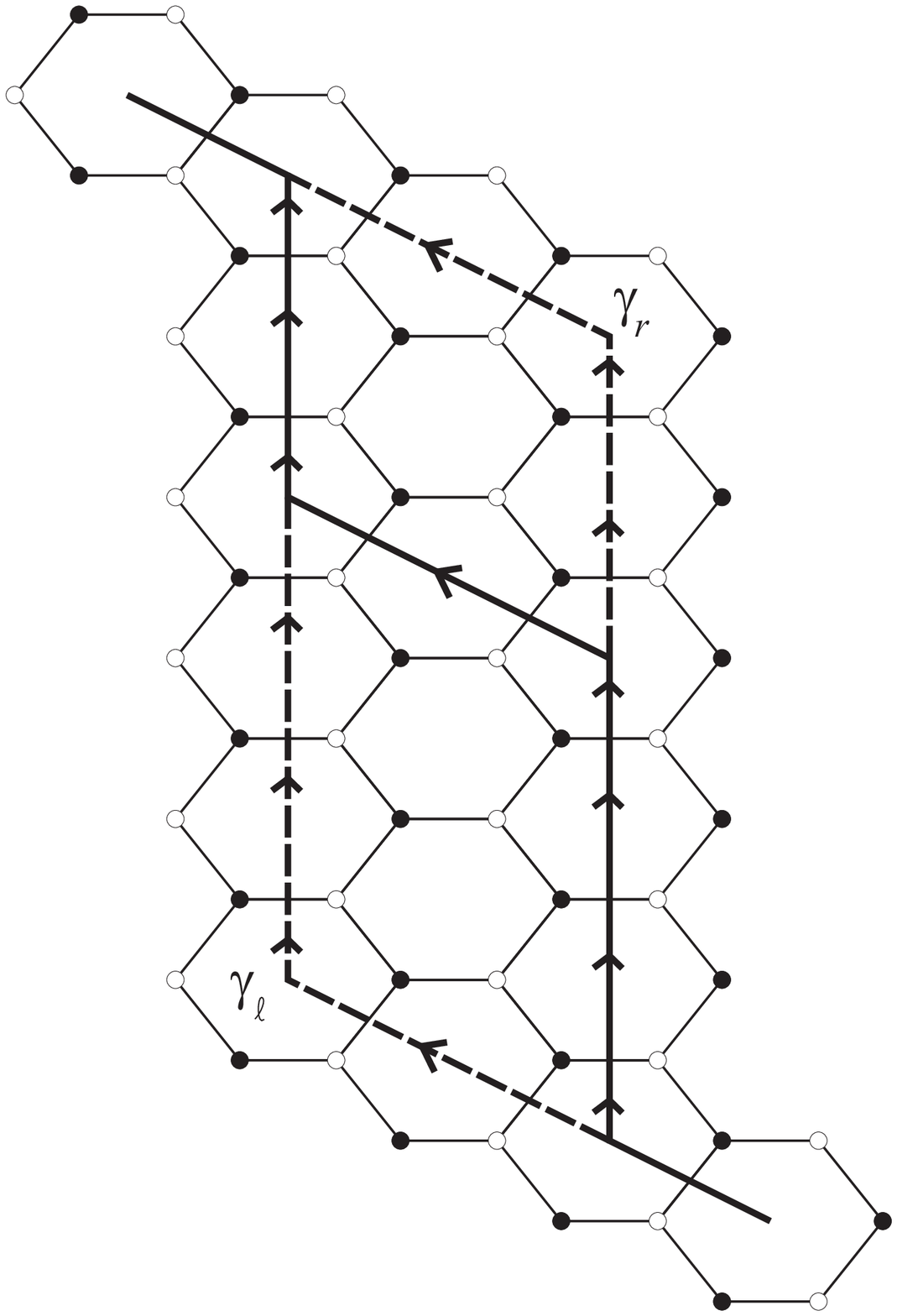}\quad\quad
    \includegraphics[scale=.33]{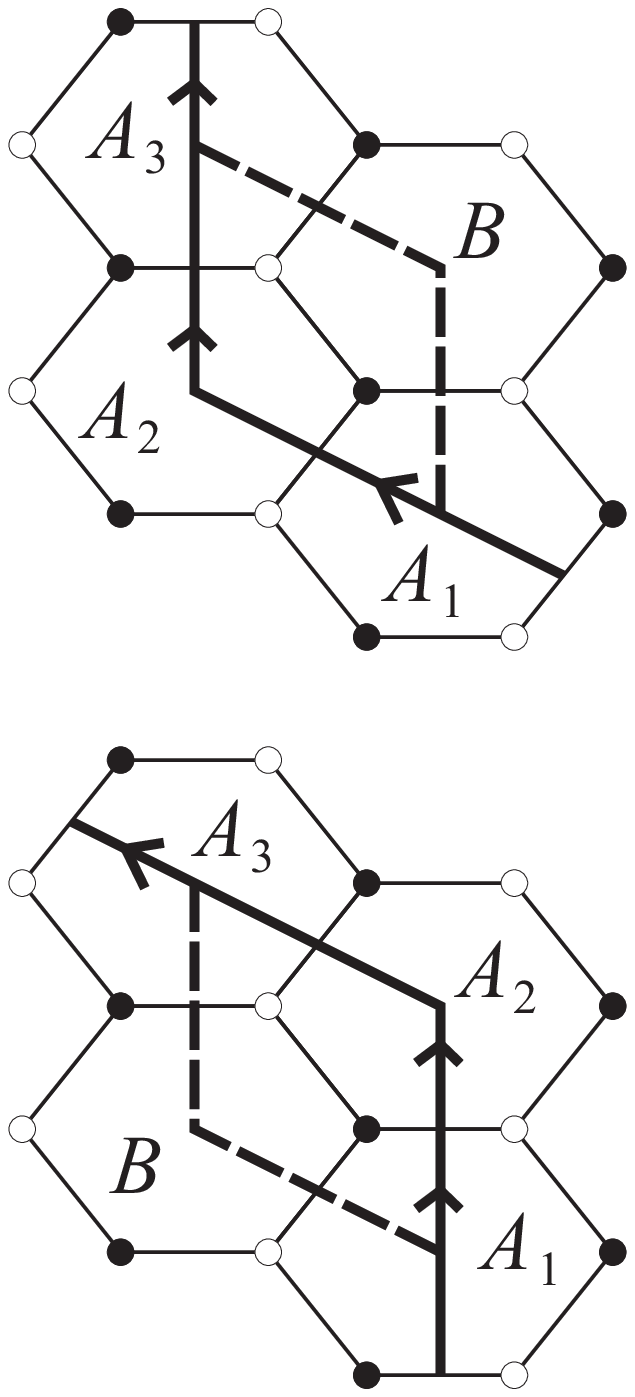}
    \end{center}
\mycap{Reducing the number of turns. \label{smallp:fig}}
\end{figure}

Our next claim is the following:
\begin{itemize}
\item[$(3)$] There exists $\gamma\in\Gamma_n$ described by a word $(H_\ell H_r)^p$ with $p\leqslant 1$.
\end{itemize}
By $(2)$ and the fact that subwords $H_\ell H_\ell$ or $H_rH_r$ are impossible in $\gamma\in\Gamma_n$,
we have a $\gamma\in\Gamma_n$ described by a word $(H_\ell H_r)^p$. Now suppose $p\geqslant 2$, consider a portion
of $\gamma$ described by $H_rH_\ell H_r H_\ell$ as in Fig.~\ref{smallp:fig}-left
and try to construct the two loops $\gamma_\ell$ and $\gamma_r$
as in Fig.~\ref{smallp:fig}-centre by repeated application of the moves in Fig.~\ref{smallp:fig}-right.
If one of $\gamma_\ell$ or $\gamma_r$ exists it belongs to $\Gamma_n$ and it is described by $(H_\ell H_r)^{p-1}$, so we can
conclude recursively. The construction of $\gamma_\ell$ or $\gamma_r$ may fail only if when we apply an elementary move as in Fig.~\ref{smallp:fig}-right to $\alpha\in\Gamma_n$ the region $B$ is\ldots
\begin{itemize}
  \item the square $S$; this would contradict (1), so it is impossible;
  \item already in $\alpha$; but then $B$ is not one of $A_1,A_2,A_3$ because all regions are embedded, and
  it easily follows that $\alpha$ is homologous to the sum of two shorter loops, which is absurd because at least one of them
  would be non-trivial;
  \item the octagon $O$; this is indeed possible, but it cannot happen both to the left and to the right, otherwise we would
  get a simplicial loop in $\Dhat$ intersecting $\gamma$ transversely at one point,
  whence non-trivial, and shorter than $\gamma$ (actually, already at least by $1$ shorter than the portion of $\gamma$
  described by $H_rH_\ell H_r H_\ell$).
\end{itemize}

\begin{figure}
    \begin{center}
    \includegraphics[scale=.33]{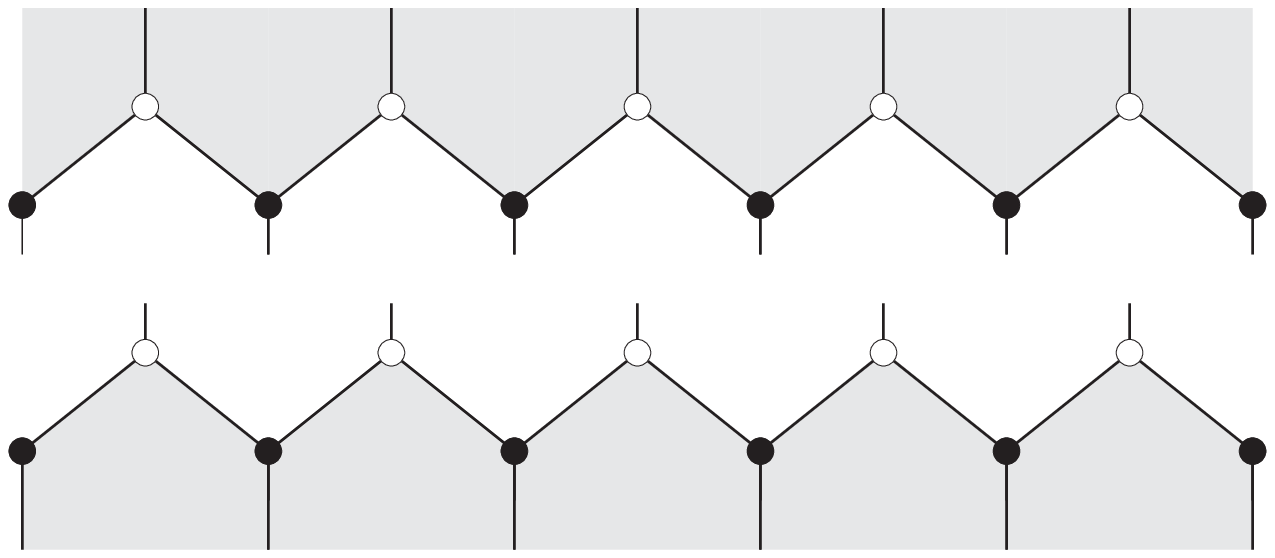}\quad\quad
    \includegraphics[scale=.33]{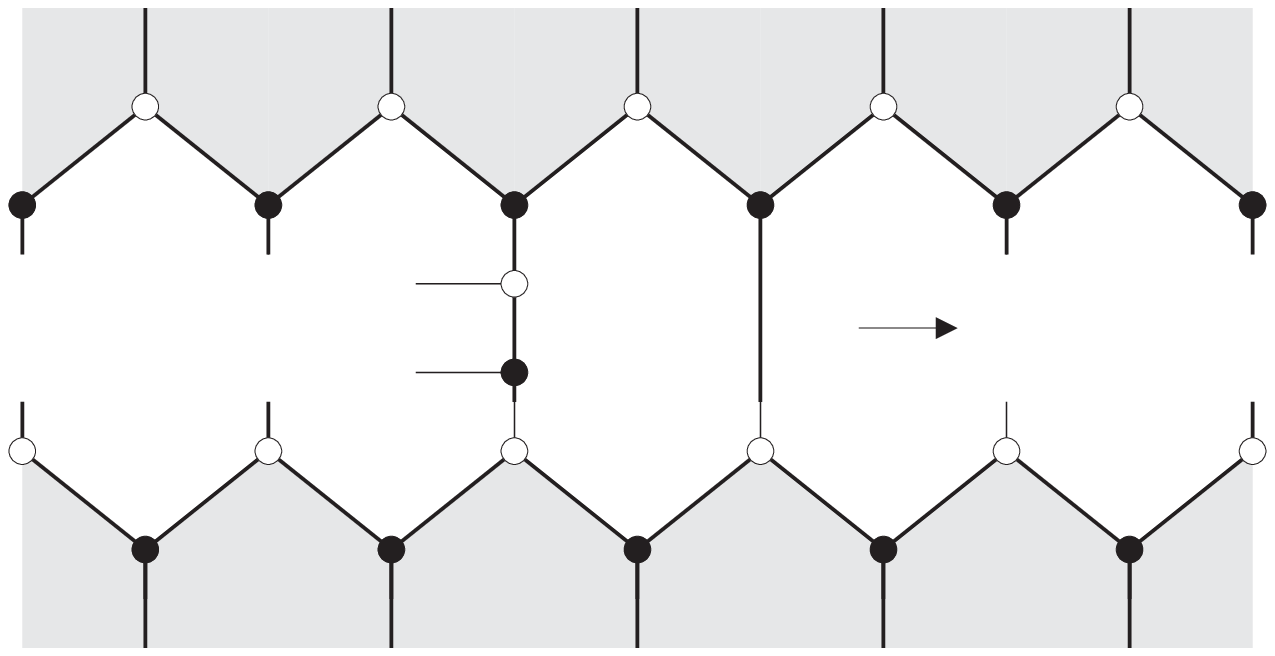}
    \end{center}
\mycap{Conclusion for the shape $H_d^q$. \label{straight:fig}}
\end{figure}

\begin{figure}
    \begin{center}
    \includegraphics[scale=.33]{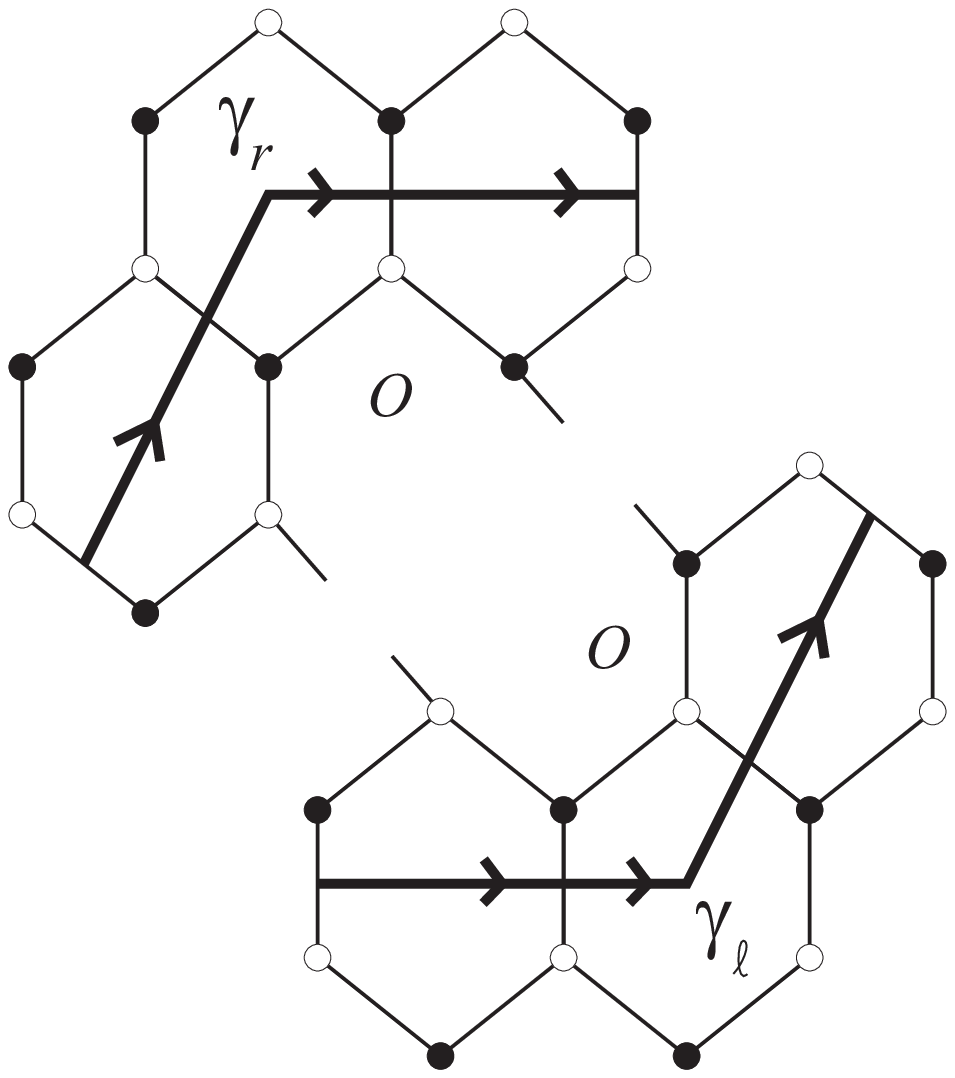}\quad\quad
    \includegraphics[scale=.33]{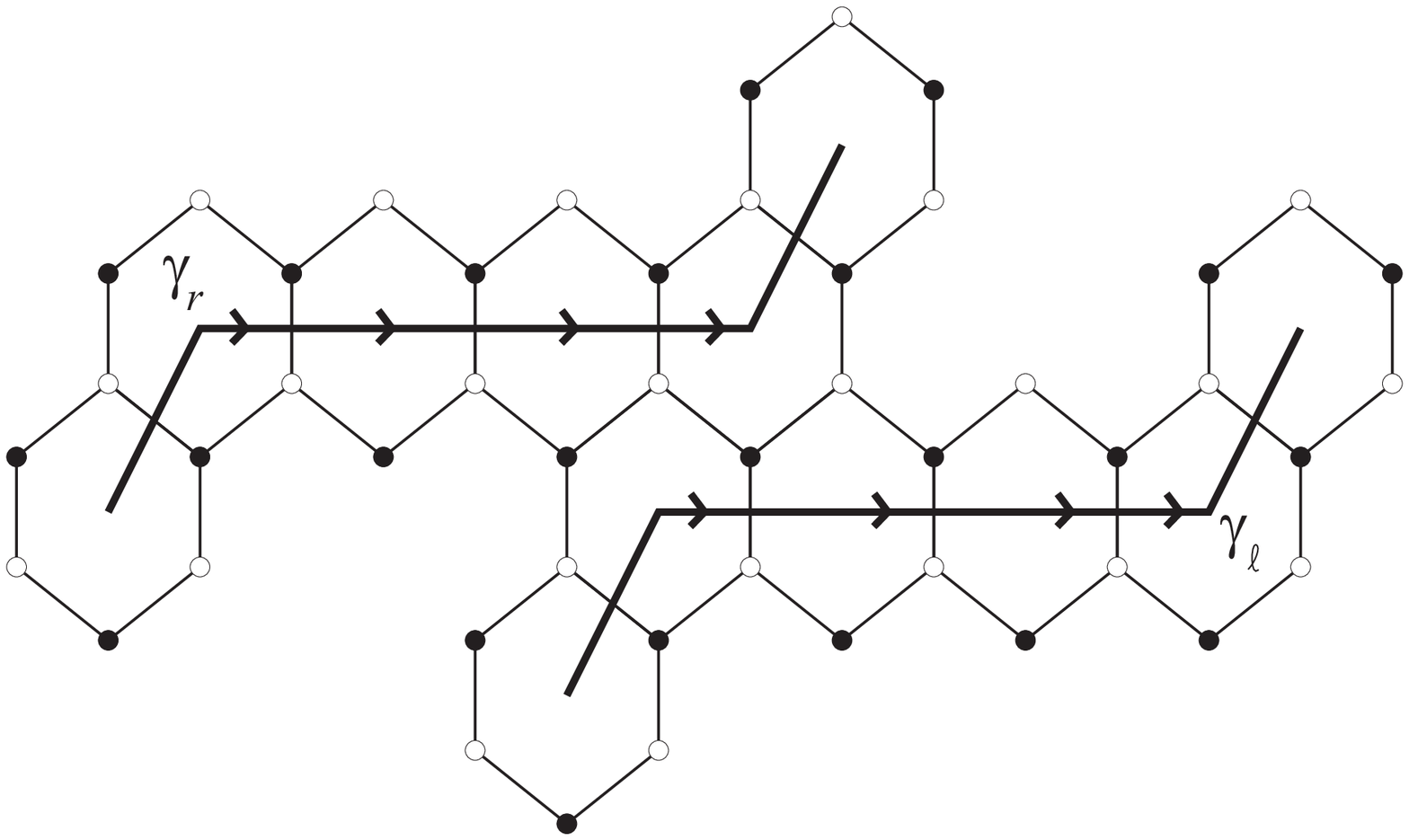}\\ \ \\
    \includegraphics[scale=.29]{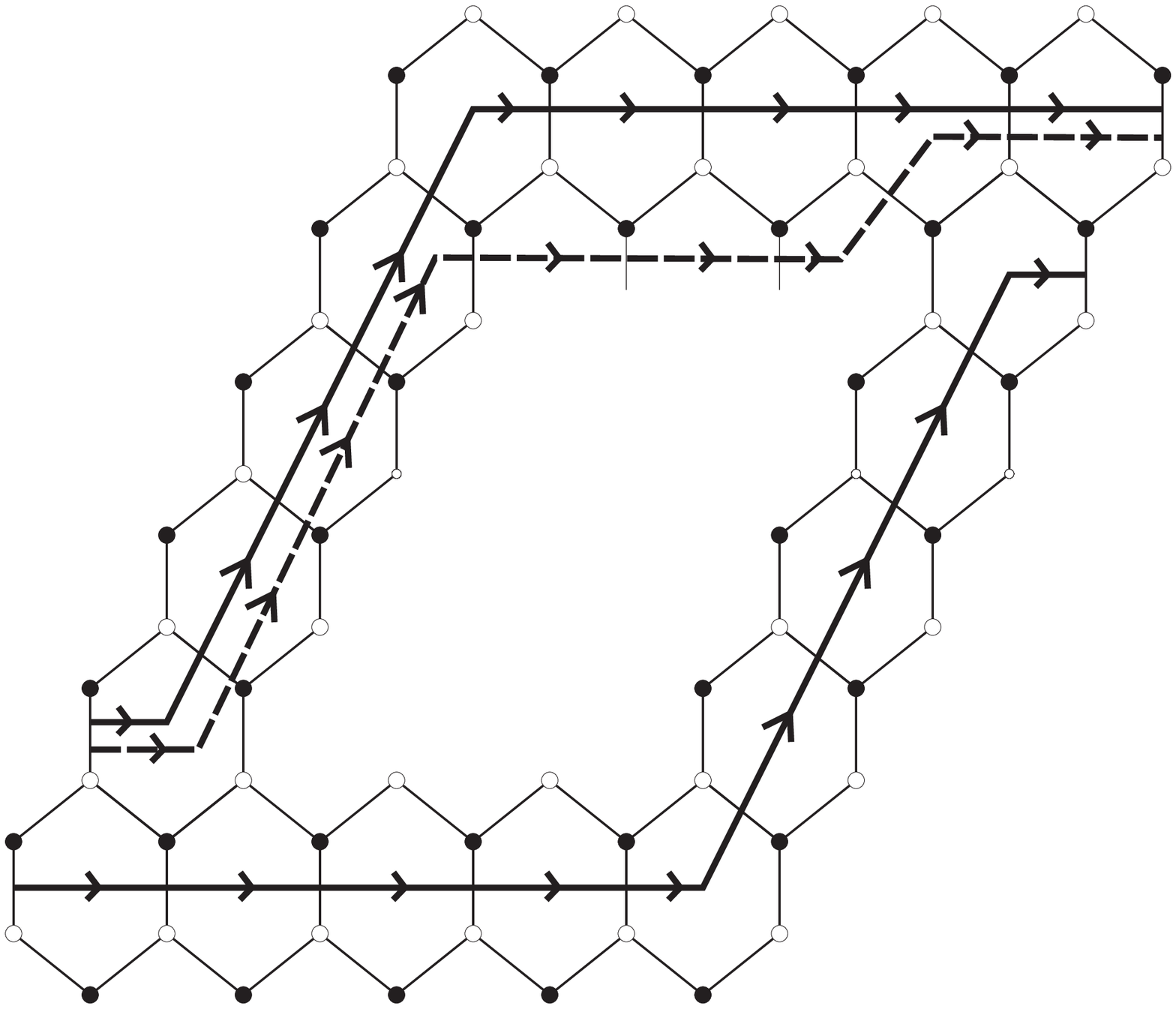}\quad\quad
    \includegraphics[scale=.29]{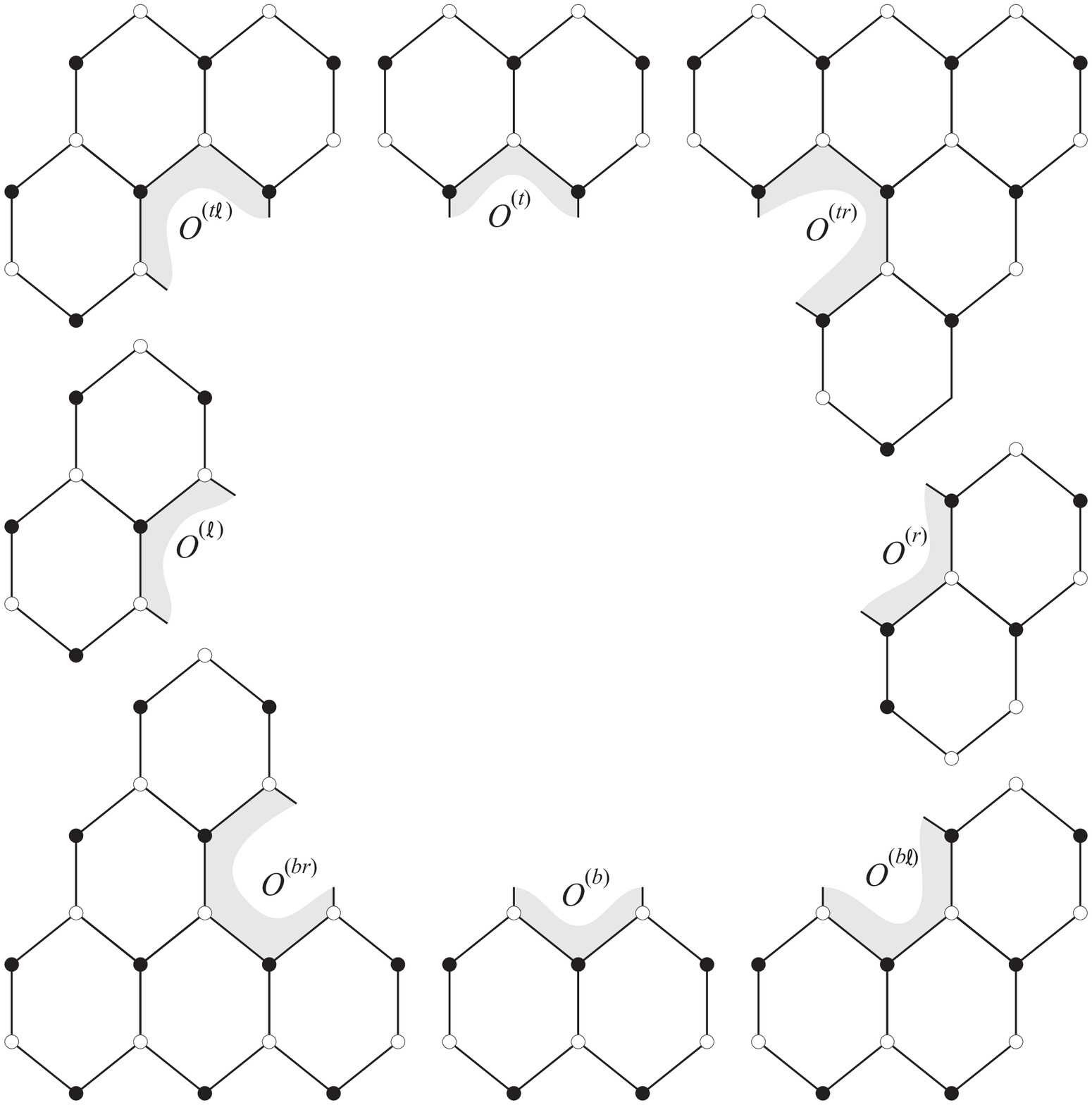}
    \end{center}
\mycap{Conclusion for the shape $H_\ell H_d^q H_r H_d^t$. \label{curved:fig}}
\end{figure}

We now include again the $H_d$'s in the notation for the word describing a loop.
It follows from (3) that there exists $\gamma\in\Gamma_n$ of shape $H_d^q$ or $H_\ell H_d^q H_r H_d^t$.
To conclude the proof we set $\gamma_\ell=\gamma_r=\gamma$ and we apply to $\gamma_\ell$ and $\gamma_r$ as long as possible
the following moves (that we describe for $\gamma_\ell$ only):
\begin{itemize}
\item If $O$ is not incident to the left margin of $\gamma_\ell$ we entirely push $\gamma_\ell$
to its left, as in Fig.~\ref{moves:fig}-left;
\begin{figure}
    \begin{center}
    \includegraphics[scale=.33]{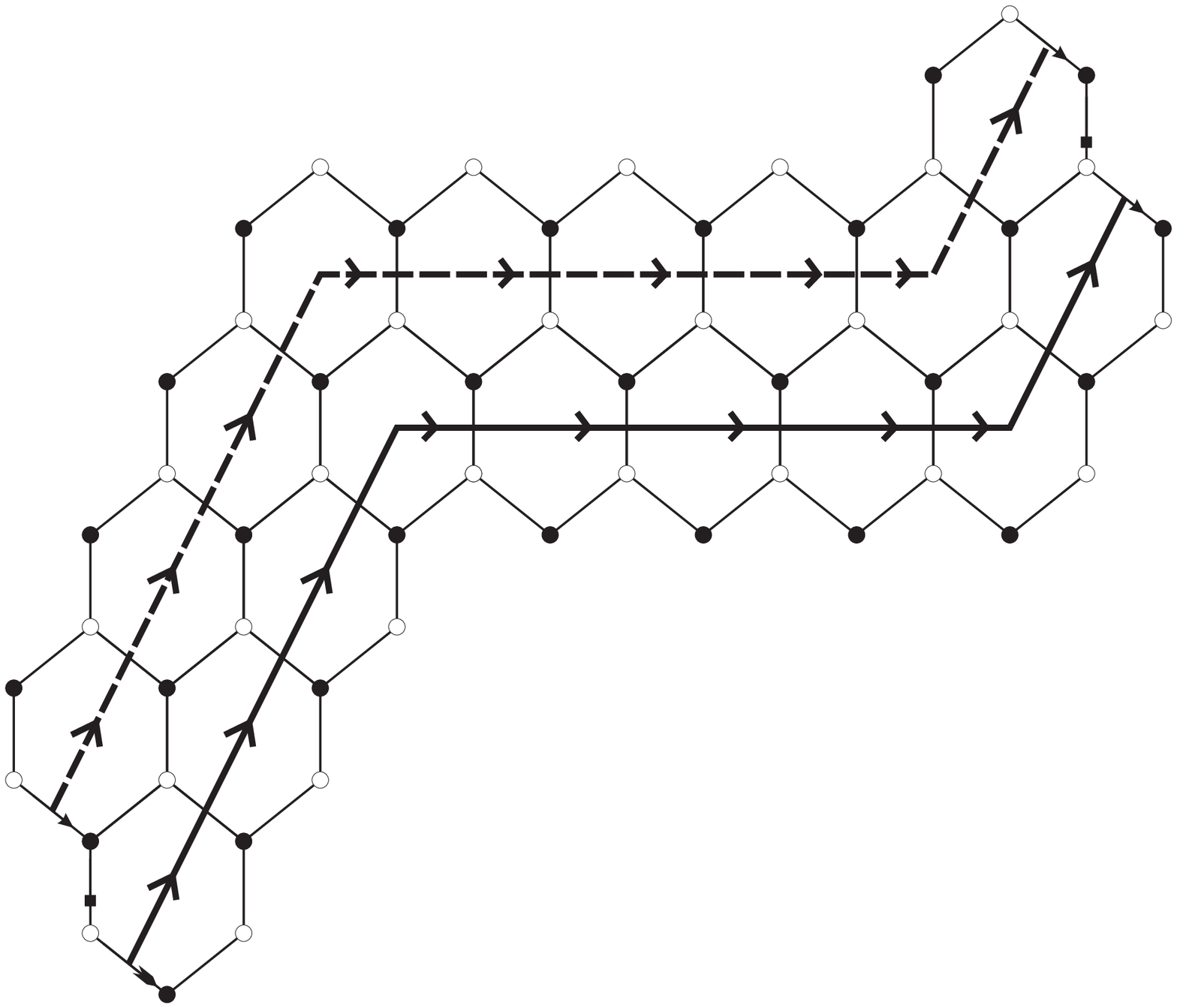}\hspace{-.5cm}
    \includegraphics[scale=.33]{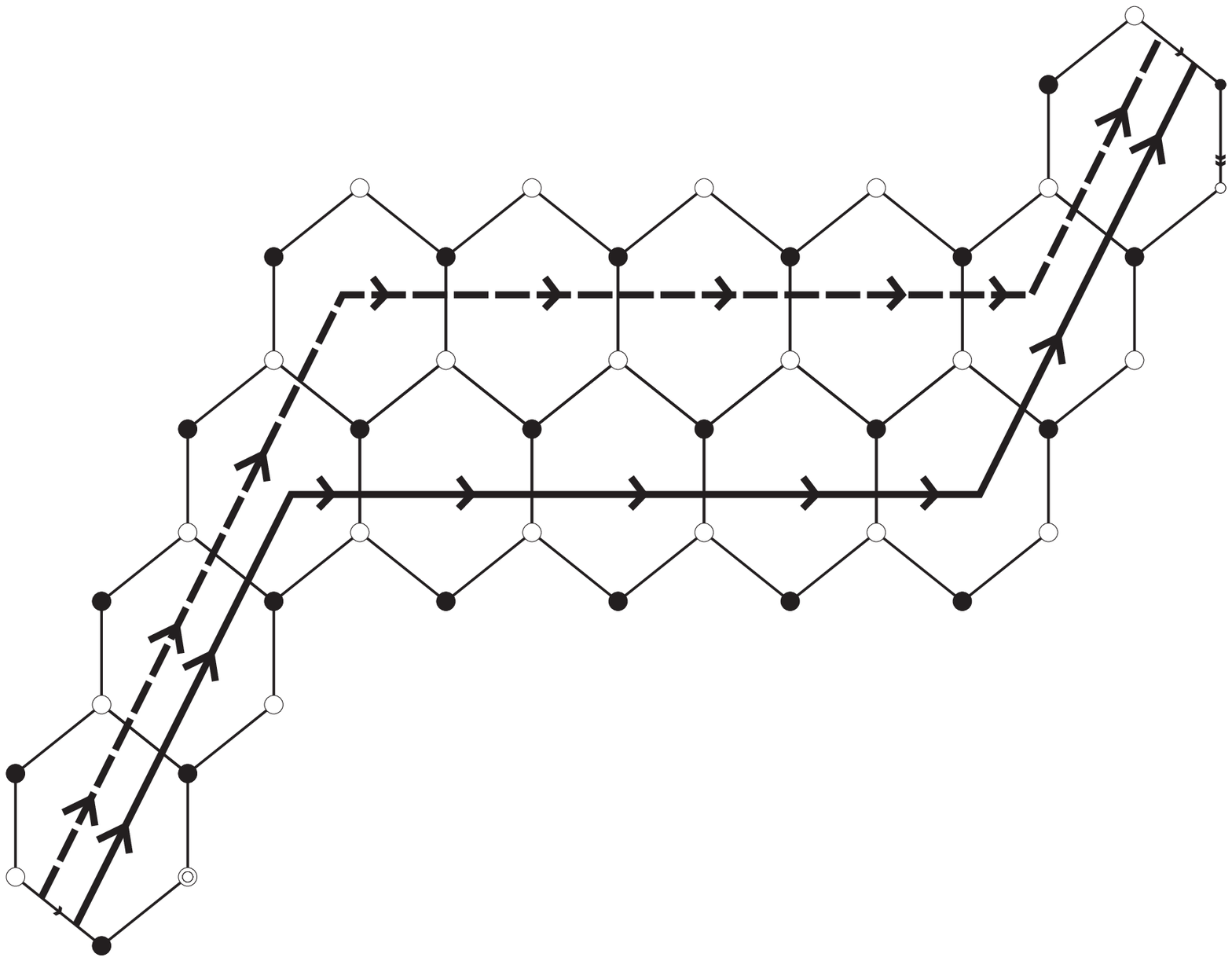}
    \end{center}
\mycap{Evolution of $\gamma_\ell$. The second move is performed so that $O$ is not included in the new loop. \label{moves:fig}}
\end{figure}
\item If $\gamma$ has shape $H_\ell H_d^q H_r H_d^t$ and $O$ is incident to the left margin of $\gamma_\ell$ but
not to $H_\ell$, we note that $O$ is not incident to either $H_d^tH_\ell$ or to $H_\ell H_d^q$, and we
partially push $\gamma_\ell$ to its left so not to include $O$, as in Fig.~\ref{moves:fig}-right (this is the case where $O$ is
not incident to $H_d^tH_\ell$).
\end{itemize}
Note that by construction the new $\gamma_\ell$ does not contain $O$, so it also does not contain $S$ by (1),
hence it has the same shape $H_d^q$ or $H_\ell H_d^q H_r H_d^t$ as the old $\gamma_\ell$. Therefore at any time
$\gamma_\ell$ and $\gamma_r$ have the same shape as the original $\gamma$. We stop applying the moves when one of the
following situations is reached:
\begin{itemize}
  \item[(a)] The left margin of $\gamma_\ell$ and the right margin of $\gamma_r$ overlap;
  \item[(b)] $\gamma_\ell$ and $\gamma_r$ have shape $H_d^q$ and $O$ is incident to
  the left margin of $\gamma_\ell$ and to the right margin of $\gamma_r$;
  \item[(c)] $\gamma_\ell$ and $\gamma_r$ have shape $H_\ell H_d^q H_r H_d^t$ and $O$ is incident to
  the left margin of $\gamma_\ell$ in $H_\ell$ and to the right margin of $\gamma_r$ in $H_r$.
\end{itemize}
Case (a) with $\gamma_\ell$ and $\gamma_r$ of shape $H_d^q$ is impossible, because the left margin
of $\gamma_\ell$ and the right margin of $\gamma_r$ would close up like a zip, leaving no space for $S$ and $O$, see Fig.~\ref{straight:fig}-left.
We postpone the treatment of case (a) with $\gamma_\ell$ and $\gamma_r$ of shape $H_\ell H_d^q H_r H_d^t$, to face the easier cases (b) and (c).
For (b), we have the situation of Fig.~\ref{straight:fig}-right, where in the direction given by the arrow we must have a strip
of identical hexagons that can never close up. Case (c), excluding (a), is trivial: the region that should be $O$ cannot close up
with fewer than $10$ vertices, see Fig.~\ref{curved:fig}-top/left.

In case (a) for the shape $H_\ell H_d^q H_r H_d^t$, the left margin of $\gamma_\ell$ can overlap with
the right margin of $\gamma_r$ only along a segment
as in Fig.~\ref{curved:fig}-top/right ---this segment has type $H_rH_d^s$ in $\gamma_\ell$ and
$H_d^sH_\ell$ in $\gamma_r$, in particular it uses $H_r$ from $\gamma_\ell$ and $H_\ell$ from $\gamma_r$, so
there is only one.  Therefore the rest of $\gamma_\ell$ and $\gamma_r$ delimit an
$x\times y$ rhombic area $R$ as in Fig.~\ref{curved:fig}-bottom/left (with $x\times y=3\times 4$ in the figure),
that must contain $S$ and $O$. Note that the $H$'s incident to $\partial R$ are pairwise distinct:
for the initial $\gamma$ the left margin cannot be incident to the right margin, otherwise a move as in
Fig.~\ref{gamma12moves:fig}-left/centre would contradict its minimality, and during the construction of $\gamma_\ell$ and $\gamma_r$
only new $H$'s are added.
If $O$ is not incident to one of the four sides of $R$ we can modify $\gamma_\ell$
or $\gamma_r$ as suggested already in Fig.~\ref{curved:fig}-bottom/left. This modification changes the
shape of $\gamma_\ell$ or $\gamma_r$, but:
\begin{itemize}
\item The modified loop is still minimal and does not contain $O$, so it does not contain $S$;
\item The area $R$ into which $O$ and $S$ are forced to lie remains a rhombus,
\item The $H$'s incident to $\partial R$ are pairwise distinct (otherwise $R$ closes up leaving no space for $O$ or $S$).
\end{itemize}
We can iterate this modification, shrinking
$R$ until $O$ is incident to all the four sides of $\partial R$. If $R$ is $1\times 1$ of course there is space in $R$ only for an $H$.
If $R$ is $1\times y$ or $x\times 1$ with $x,y\geqslant 2$, the fact that the $H$'s incident to $\partial R$ are distinct
implies that the vertices of $\partial R$ are distinct, so a region incident to all the four sides or $\partial R$ must have
at least $10$ vertices. If $R$ is $x\times y$ with $x,y\geqslant 2$ then $O$ contains some of the germs of regions $O^{(*)}$ in
Fig.~\ref{curved:fig}-bottom/right so as to touch all the $r/t/\ell/b$ sides of $\partial R$.
An easy analysis shows that any identification between two vertices of the $O^{(*)}$'s would force two
$H$'s incident to $\partial R$ to coincide, so it is impossible. This implies that any $O^{(*)}$
actually contained in $O$ contributes to the number of vertices of $O$ with as many vertices as one sees in
Fig.~\ref{curved:fig}-bottom/right, namely $3$ for $O^{(r)}$, $O^{(t)}$, $O^{(\ell)}$, $O^{(b)}$, then
$4$ for $O^{(t\ell)}$, $O^{(br)}$, and finally $5$ for $O^{(tr)}$, $O^{(b\ell)}$.
Therefore, a region can touch all of $r/t/\ell/b$ with a total of no more that $8$ vertices only if it includes
$O^{(t\ell)}$ and $O^{(br)}$, but then the vertex colors again imply that the number of vertices is at least $10$.
This gives the final contradiction and concludes the proof.

\noindent
\'Ecole Normale Sup\'erieure de Rennes\\
Campus de Ker Lann\\
Avenue Robert Schuman\\
35170 BRUZ -- France\\
\texttt{Tom.Ferragut@ens-rennes.fr}

\bigskip

\noindent
Dipartimento di Matematica\\
Universit\`a di Pisa\\
Largo Bruno Pontecorvo, 5\\
56127 PISA -- Italy\\
\texttt{petronio@dm.unipi.it}

\end{document}